\documentstyle[amssymb,12pt]{amsart}%\input amssym.def
\newtheorem{prop}{Proposition}[section]
\newtheorem{prop:def}{Proposition-Definition}[section]

\newtheorem{thm}{Theorem}[section]
\newtheorem{cor}{Corollary}[section]

\theoremstyle{remark}
\newtheorem{remark}{Remark}

\begin{document}
\newcommand{\nc}{\newcommand} \nc{\on}{\operatorname}
\nc{\pa}{\partial}
\nc{\cA}{{\cal A}}\nc{\cB}{{\cal B}}\nc{\cC}{{\cal C}}
\nc{\cE}{{\cal E}}\nc{\cG}{{\cal G}}\nc{\cH}{{\cal H}}
\nc{\cX}{{\cal X}}\nc{\cR}{{\cal R}}\nc{\cL}{{\cal L}}
\nc{\cK}{{\cal K}}
\nc{\sh}{\on{sh}}\nc{\Id}{\on{Id}}\nc{\Diff}{\on{Diff}}
\nc{\ad}{\on{ad}}\nc{\Der}{\on{Der}}\nc{\End}{\on{End}}
\nc{\res}{\on{res}}\nc{\ddiv}{\on{div}}
\nc{\card}{\on{card}}\nc{\dimm}{\on{dim}}
\nc{\Jac}{\on{Jac}}\nc{\Ker}{\on{Ker}}
\nc{\Imm}{\on{Im}}\nc{\limm}{\on{lim}}\nc{\Ad}{\on{Ad}}
\nc{\ev}{\on{ev}}
\nc{\Hol}{\on{Hol}}\nc{\Det}{\on{Det}}
\nc{\de}{\delta}\nc{\si}{\sigma}\nc{\ve}{\varepsilon}
\nc{\al}{\alpha}\nc{\vp}{\varphi}
\nc{\CC}{{\mathbb C}}\nc{\ZZ}{{\mathbb Z}}
\nc{\NN}{{\mathbb N}}\nc{\zz}{{\mathbf z}}\nc{\nn}{{\mathbf n}}
\nc{\AAA}{{\mathbb A}}\nc{\cO}{{\cal O}} \nc{\cF}{{\cal F}}
\nc{\la}{{\lambda}}\nc{\G}{{\mathfrak g}}\nc{\mm}{{\mathbf m}}
\nc{\A}{{\mathfrak a}}
\nc{\HH}{{\mathfrak h}}
\nc{\N}{{\mathfrak n}}\nc{\B}{{\mathfrak b}}
\nc{\La}{\Lambda}
\nc{\g}{\gamma}\nc{\eps}{\epsilon}\nc{\wt}{\widetilde}
\nc{\wh}{\widehat}
\nc{\bn}{\begin{equation}}\nc{\en}{\end{equation}}
\nc{\SL}{{\mathfrak{sl}}}\nc{\ttt}{{\mathfrak{t}}}

% ****** GISPIC **********
%
%** by GISLI MASON *******
%
%**for commutative diagrams
%

\newcommand{\ldar}[1]{\begin{picture}(10,50)(-5,-25)
\put(0,25){\vector(0,-1){50}}
\put(5,0){\mbox{$#1$}} 
\end{picture}}

\newcommand{\lrar}[1]{\begin{picture}(50,10)(-25,-5)
\put(-25,0){\vector(1,0){50}}
\put(0,5){\makebox(0,0)[b]{\mbox{$#1$}}}
\end{picture}}

\newcommand{\luar}[1]{\begin{picture}(10,50)(-5,-25)
\put(0,-25){\vector(0,1){50}}
\put(5,0){\mbox{$#1$}}
\end{picture}}

\title[Correlation functions of Drinfeld currents and shuffle
algebras] {On correlation functions of Drinfeld currents ans shuffle
  algebras}

\author{B. Enriquez}

\address{Centre de Math\'ematiques, Ecole Polytechnique, UMR 7640 du
  CNRS, 91128 Palaiseau, France}

\address{D-Math, FIM, ETH-Zentrum, HG G46, CH-8092 Z\"urich,
  Switzerland}

%\author{V. Rubtsov}

%\address{V.R.: D\'ept. de Math\'ematiques, Univ. d'Angers, 
%2, Bd. Lavoisier, 49045 Angers, France}

%\address{ITEP, 25 Bol. Cheremushkinskaia, 117259 Moscou, Russie}

\date{September 1998}

\begin{abstract}
  We express the vanishing conditions satisfied by the correlation
  functions of Drinfeld currents of quantum affine algebras, imposed
  by the quantum Serre relations. We discuss the relation of these
  vanishing conditions with a shuffle algebra description of the
  algebra of Drinfeld currents.
\end{abstract}

\maketitle

\subsection*{Introduction}
This paper is concerned with the functional properties of the
correlation functions of Drinfeld currents of untwisted quantum affine
algebras.  Let $\bar\G$ be a simple finite-dimensional complex Lie
algebra, $\G$ be its affinization and $U_q\G$ the associated quantum
algebra. We will denote the positive nilpotent currents of $U_q\G$ by
$e_\al(z)$, where $z$ is a formal variable and $\al$ is a simple root
of $\bar\G$.  The correlation functions of these currents are defined
as follows.

For $V$ a highest weight module over $U_q\G$, $v$ in $V$ and $\xi$ in
$V^*$ weight vectors, one considers the series
\begin{equation} \label{def:corr} 
  f(z^{(\al)}_i) = \langle \xi, \prod_{\al = 1}^r \prod_{i\in
    I_\al} e_\al(z_i^{(\al)})v\rangle,
\end{equation} 
where the $I_\al$ are the sets $\{1,\ldots,n_\al\}$, and $r$ is the
rank of $\bar\G$; this series is defined in
$\CC[[z_i^{(\al)},z_i^{(\al)-1}]]$.

In the classical case ($q=1$), the functional properties of these
correlation functions are the following. 

\begin{thm} (\cite{FS}) 
  Denote by $(a_{\al\beta})_{1\leq\al,\beta\leq r}$ the Cartan matrix
  of $\bar\G$. $f(z^{(\al)}_i)$ belongs to
  $\CC((z^{(1)}_1))((z^{(1)}_2)) \cdots ((z^{(r)}_{n_r}))$; it is the
  expansion of a rational function in the $(z^{(\al)}_i)_{i\in I_\al}$
  for each $\al$, regular except for simple poles at the diagonals
  $z^{(\al)}_i = z^{(\beta)}_j$ for $a_{\al\beta}\neq 0$, and when
  some $z^{(\al)}_i$ meets the origin; it satisfies the vanishing
  conditions
\begin{equation} \label{class:rels}
  \res_{z^{(\al)}_{i_N} = z^{(\beta)}_j} \cdots \res_{z^{(\al)}_{i_1}
    = z^{(\beta)}_j} f(z^{(\al)}_i) dz^{(\al)}_{i_1}\cdots
  dz^{(\al)}_{i_N} = 0 ,
\end{equation}
for any simple roots $\al,\beta$, where $N = 1-a_{\al\beta}$, $j$ is
some index of $I_{\beta}$ and $i_1, \ldots, i_N$ are $N$ distinct
indices of $I_\al$. In other words, $f(z^{(\al)}_i)$ has the form
$$ f(z^{(\al)}_i) = {1\over{\prod_{\al < \beta}\prod_{i\in I_\al, j\in
      I_\beta}(z^{(\al)}_i - z^{(\beta)}_j )}} A(z^{(\al)}_i),
$$ with $A(z^{(\al)}_i)$ in $\CC[z^{(\al)}_i,z^{(\al)-1}_i]$ zero
whenever $z^{(\al)}_{i} = z^{(\beta)}_{i_1} = \cdots =
z^{(\al)}_{i_{1-a_{\al\beta}}}$.
\end{thm}

The vanishing conditions are consequences of the Serre relations, and
of the identity $\langle \xi, [x,y](z) v\rangle = \res_{z'=z}\langle
\xi, x(z')y(z) v \rangle dz$, for $x,y$ in $\bar\G$ and $x(z)$ the
field associated with $x$.

Our goal in this paper is to express vanishing conditions for
correlators of Drinfeld currents, analogous to the relations
(\ref{class:rels}). We show
\begin{thm} \label{main} 
  Assume that $\bar\G$ is not of type $G_2$.  Let $q$ be an arbitrary nonzero
  complex number.  Let $d_\al$ be the symmetrizing factors of $\bar
  \G$, so that $(d_\al a_{\al\beta})$ is a symmetric matrix. Set
  $q_\al = q^{d_\al}$.  The correlation function $f(z^{(\al)}_i)$
  defined by (\ref{def:corr}) has the following properties:

  1) it belongs to $\CC((z^{(1)}_1))((z^{(1)}_2)) \cdots
  ((z^{(r)}_{n_r}))$; 

  2) it is the expansion of a rational function on
  $\prod_{\al}(\CC^{\times})^{I_\al}$, regular everywhere except for
  simple poles on shifted diagonals $q_{\al}^{a_{\al\beta}}z^{(\al)}_i
  = z^{(\beta)}_j$ for $\al < \beta$ or $\al = \beta$ and $i<j$, or
  when some $z^{(\al)}_i$ meets the origin; it satisfies the twisted
  symmetry relations
  \begin{equation} \label{tw:symm} 
    (q_\al^{2} z^{(\al)}_i- z^{(\al)}_j) f(z^{(\al)}_i) =
    (z^{(\al)}_i- q^{2}_{\al} z^{(\al)}_j)
    (\sigma^{(\al)}_{ij}f)(z^{(\al)}_i),
  \end{equation} where $\sigma^{(\al)}_{ij}$ is the interchange of variables
  $z^{(\al)}_i$ and $z^{(\al)}_j$. In other words, it has the form
  \begin{equation} \label{can} 
    f(z^{(\al)}_i) = {  { \prod_\al \prod_{i,j\in I_\al,
        i<j}(z^{(\al)}_i - z^{(\al)}_j )}\over{ \prod_ {\al\leq\beta}
        \prod_{i\in I_\al,j\in I_\beta, i< j \on{if} \al = \beta}
        (q_{\al}^{a_{\al\beta}} z^{(\al)}_i - z^{(\beta)}_j})}
    A(z^{(\al)}_i) , \end{equation} where $A$ belongs to the space of
  Laurent polynomials $\CC[z^{(\al)}_i,z^{(\al)-1}_i]$ and is
  symmetric in the $(z^{(\al)}_i)_{i\in I_\al}$, for each $\al$;

  3) moreover, $A$ satisfies
  \begin{equation} \label{vanishing} A(z^{(\al)}_i) = 0 \quad 
    \on{when} \quad z^{(\al)}_{i_1} = q_{\al}^2 z^{(\al)}_{i_2} =
    \ldots = q_{\al}^{-2a_{\al\beta}}z^{(\al)}_{i_{1-a_{\al\beta}}}
    = q_{\al}^{-a_{\al\beta}} z^{(\beta)}_j, 
  \end{equation}
  for any $\al,\beta$, where $j$ belongs to $I_\beta$ and the $i_j$
  are pairwise different elements of $I_\al$. In other words, $f$ has
  the form
  $$ f(z^{(\al)}_i) = {{ \prod_\al \prod_{i,j\in I_\al,
        i<j}(z^{(\al)}_i - z^{(\al)}_j )}\over{ \prod_{\al\leq \beta,
        a_{\al\beta}\neq 0} \prod_{i\in I_\al,j\in I_\beta, i<
        j\on{if}\al = \beta} (q_\al^{a_{\al\beta}}z^{(\al)}_i -
      z^{(\beta)}_j}) } B(z^{(\al)}_i),
  $$ where $B(z^{(\al)}_i)$ is in $\CC[z^{(\al)}_i, z^{(\al)-1}_i]$
  and satisfies (\ref{vanishing}) for all $\al,\beta$ such that
  $a_{\al\beta}\neq 0$.
\end{thm}

1) and 2) are standard facts that are explained in sect.
\ref{sect:std}. The proof of 3) (sect. \ref{sect:proof:3}) rests on
some delta-function identities that are established in sect.
\ref{sect:comb:ids}. 

We then discuss the relation of this result with the shuffle algebra
description of the algebra $U_q\N_+$ (sect.  \ref{sect:fun:desc}).
Define $\overline{Sh}$ the direct sum
$\oplus_{\nn\in\NN^r}\overline{Sh}_\nn$, where for $\nn =
(n_1,\ldots,n_r)$, $\overline{Sh}_\nn$ is the subspace of
$\CC((z^{(1)}_1))\cdots ((z^{(r)}_{n_r}))$ formed of the elements of
the form (\ref{can}), where $A$ belongs to
$\CC[z^{(\al)}_i,z^{(\al)-1}_i]$ and is symmetric in the group of
variables $(z^{(\al)}_i)_{i\in I_\al}$, for each $\al$.
$\overline{Sh}$ is endowed with the a product of shuffle type, which
was introduced by B. Feigin and A. Odesskii (\cite{FO}). We show

\begin{thm} \label{thm:shuffle}
  There exists an algebra morphism $\pi$ from $U_q\N_+$ to
  $\overline{Sh}$, sending each $e_\al[n]$ to $(z^{(\al)})^n\in
  \overline{Sh}_{\delta_\al}$ (here $\delta_\al$ is the vector with
  coordinates $(\delta_\al)_\beta = \delta_{\al\beta}$).

  If $\bar\G$ is not of type $G_2$, $\pi$ maps $U_q\N_+$ to
  the subspace $Sh$ of $\overline{Sh}$, defined as the direct sum
  $\oplus_{\nn\in\NN^r}Sh_\nn$, where $Sh_\nn$ consists of the
  elements of the form (\ref{can}), such that $A$ satisfies the
  vanishing conditions (\ref{vanishing}).
\end{thm}

It is natural to expect that $\pi$ is an isomorphism, so that $Sh$ is
actually a subalgebra of $\overline{Sh}$. This result would probably
lead to a simple proof of the PBW theorem for $U_q\N_+$ (see
\cite{Beck,ER}).

We thank E. Frenkel, who suggested that Thm. \ref{main} has some
connection with integral formulas for the qKZ equations in higher
rank. We would also like to thank N. Andruskiewitsch and G. Felder for
discussions related to this paper, and the FIM, ETHZ, for hospitality
during the preparation of this work.

\section{Properties 1) and 2) of the correlation functions}
\label{sect:std}

Recall first the presentation of the positive nilpotent part $U_q\N_+$
of $U_q\G$.  It has generators $e_\al[n]$, $\al$ simple, $n$ integer,
organized in generating series $e_\al(z) = \sum_{n\in \ZZ}
e_\al[n]z^{-n}$, subject to the vertex relations
\begin{equation} \label{vertex}
  (q_\al^{a_{\al\beta}}z - w) e_\al(z)e_\beta(w) = (z -
  q_\al^{a_{\al\beta}}) e_\beta(w) e_\al(z),
\end{equation}
and to the the quantum Serre relations 
\begin{align} \label{q:Serre}
  \sum_{r=0}^{1-a_{\al\beta}} (-1)^r \bmatrix 1 - a_{\al\beta} \\ r
  \endbmatrix_{q_\al} \on{Sym}_{z_1\cdots z_{1-a_{\al\beta}}}
  (e_\al(z_1)\cdots e_\al(z_r) e_\beta(w) e_\al(z_{r+1})\cdots &
  e_\al(z_{1-a_{\al\beta}})) \\ & \nonumber = 0,
\end{align}
where $\bmatrix n \\ p \endbmatrix_q =
{{[n]^!_q}\over{[p]^!_q[n-p]^!_q}}$, and $[n]^!_q = [1]_q[2]_q\cdots
[n]_q$, with $[n]_q = {{q^n - q^{-n}}\over{q-q^{-1}}}$.

Let us pass to the proof of properties 1) and 2) of the
$f(z^{(\al)}_i)$. 

Recall that $\CC[[z_1,z_1^{-1},\cdots,z_n,z_n^{-1}]]$ is defined as
the space of expansions 
$$\sum_{i_1\ldots i_n}u_{i_1\ldots
  i_n}z_1^{i_1}\cdots z_n^{i_n},
$$ with no restrictions on the $(i_j)$.  On the other hand, for any
ring $R$, $R((z))$ is the ring defined as the localization
$R[[z]][z^{-1}]$, that is the set of formal series $\sum_{i\in\ZZ} r_i
z^{i}$, where the $r_i$ belong to $R$ and vanish if $i$ is smaller
than some integer.

For $x_i$ a sequence of generating fields, the series 
\begin{equation} \label{ex:corr}
  \langle \psi, x_1(z_1) \ldots x_n(z_n) v \rangle
\end{equation} 
is defined as an element of $\CC[[z_1,z_1^{-1},\cdots,z_n,z_n^{-1}]]$. 

Let $V$ be a highest weight module over $U_q\G$. $V$ has a direct sum
decomposition $V = \sum_{i\leq i_0} V_i$, which makes it a graded
module over $U_q \G$, endowed with the homogeneous gradation (where
each $e_\al[n]$ has homogeneous degree $n$).

(\ref{ex:corr}) is expanded as $\sum_{i_1,\ldots,\i_n\in\ZZ} \langle
\xi, x_1[i_1] \cdots x_n[i_n] v\rangle z_1^{-i_1}\cdots z_n^{-i_n}$.
On the other hand, $v$ belongs to some $V_j$. The coefficient of
$z_1^{-i_1}\cdots z_n^{-i_n}$ therefore vanishes if $i_n >N-j$ or if
$i_{n-1} > N-(j+i_n)$,... or if $i_{p}> N - (j + i_n + \cdots +
i_{p+1})$. This means that (\ref{ex:corr}) belongs to
$\CC((z_1))\cdots ((z_n))$, proving 1) of Thm. \ref{main}.

Let now show first that it implies that $f(z^{(\al)}_i)$ has the form
(\ref{can}), with $A$ in $\CC[[z^{(\al)}_i]][z^{(\al)-1}_i]$.  Rename
the fields and variables occurring in the definition of
$f(z^{(\al)}_i)$ as 
$$
f(z_1,\cdots,z_N) =\langle \xi, x_1(z_1)\cdots x_N(z_N) v\rangle. 
$$ We have the relations $(q^{b_{ij}}z-w)x_i(z)x_j(w) =
(z-q^{b_{ij}}w)x_j(w)x_i(z)$. Set for $\sigma$ in $S_n$,
$$
f_\sigma(z_1,\cdots,z_N) = \langle \xi, x_{\si(1)}(z_{\si(1)})
\cdots x_{\si(N)}(z_{\si(N)}) v \rangle. 
$$ Then $f_{\si}(z_1,\cdots,z_N)$ belongs to $\CC((z_{\si(1)})) \cdots
((z_{\si(N)}))$. Since
$$ \prod_{i<j}(q^{b_{ij}}z_i - z_j) f(z_1,\cdots,z_N) =
\prod_{\si(i)>\si(j)}(z_i -
q^{b_{ij}}z_j)\prod_{\si(i)<\si(j)}(q^{b_{ij}}z_i - z_j)
f_{\si}(z_1,\cdots,z_N), 
$$ $\prod_{i<j}(q^{b_{ij}}z_i - z_j) f(z_1,\cdots,z_N)$ belongs to the
intersection of all the $\CC((z_{\si(1)}))\cdots((z_{\si(N)}))$, with
$\si$ in $S_n$, which is $\CC[[z_1\ldots z_N]][z_1^{-1}\ldots
z_N^{-1}]$. Therefore $f_\sigma(z_1,\cdots,z_N)$ has the form
$$ f_\sigma(z_1,\cdots,z_N) =
{{B(z_1,\ldots,z_N)}\over{\prod_{i<j}(q^{b_{ij}}z_i - z_j)}}, 
$$ with $B(z_1,\ldots,z_N)$ in $\CC[[z_1\ldots z_N]][z_1^{-1}\ldots
z_N^{-1}]$ (expansion for $z_N << z_{N-1}<< \ldots$).

Now $\prod_{\al\leq\beta}\prod_{i\in I_\al,j\in I_\beta, i< j
  \on{if}\al = \beta} (q_\al^{a_{\al\beta}} z^{(\al)}_i -
z^{(\beta)}_j) f(z^{(\al)}_i)$ is totally antisymmetric in each group
of variables $(z^{(\al)}_i)_{i\in I_\al}$, for each $\al$, which
implies that $B$ has the form $\prod_{\al}\prod_{i<j}(z^{(\al)}_i -
z^{(\al)}_j)A(z^{(\al)}_i)$, with $A(z^{(\al)}_i)$ in
$\CC[z^{(\al)}_i][[z^{(\al)-1}_i]]$.

Then the fact that $\xi$ has fixed homogeneous degree implies that
$f(z^{(\al)}_i)$ has a fixed total degree in the variables
$z^{(\al)}_i$. It follows that $A$ has also a fixed total degree in
these variables. Write $A(z^{(\al)}_i) = \sum_{i^{(\al)}_j}
A[i^{(\al)}_j](z^{(\al)}_j)^{i^{(\al)}_j}$, the $A[i^{(\al)}_j]$
therefore vanish unless the sum of all $i^{(\al)}_j$ is equal to a fixed
number; they also vanish unless the $i^{(\al)}_j$ are greater than a
fixed number, which implies that all but finitely many
$A[i^{(\al)}_j]$ vanish. This means that $A$ belong to
$\CC[z^{(\al)}_i, z^{(\al)-1}_i]$, therefore proving Thm. \ref{main},
2).

\section{Delta-function identities} \label{sect:comb:ids}

The proof of Thm. \ref{main}, 3), relies on the following
combinatorial identities. Here and below, we use the convention that
${1\over{z-w}} = \sum_{i\geq 0}w^i z^{-i-1}$ and define $\delta(z,w)$
as $\sum_{i\in \ZZ}z^i w^{-i-1}$. We have $\delta(z,w) = {1\over{z-w}}
+ {1\over{w-z}}$.

\begin{prop} \label{delta}
We have for $m=1,2$ the identity
\begin{align*} \label{conj}
  & \on{Sym}_{z_1\ldots z_{m+1}} \sum_{k=0}^m \bmatrix m+1 \\ k
  \endbmatrix_q \\ & {1\over{(q^{-m}z_1 - w)\cdots (q^{-m}z_k -
      w)(q^{-m}w-z_{k+1})\cdots (q^{-m}w-z_{m+1})}} {{\prod_{i<j}(z_i
      - z_j)}\over{\prod_{i<j}(q^2 z_i - z_j)}} \\ & = q^m
  \on{Sym}_{z_1\cdots z_{m+1}} \delta(w,q^{-m}z_1)\delta(z_1,q^2 z_2)
  \cdots \delta(z_m,q^2 z_{m+1}), 
\end{align*}
where we set 
$$ \on{Sym}_{z_1\cdots z_{m+1}} f(z_1,\ldots,z_{m+1}) = \sum_{\si\in
  S_{m+1}} f(z_{\si(1)},\ldots,z_{\si(m+1)}).
$$
\end{prop}

{\em Proof of Prop. \ref{delta}.} In the case $m=1$, the left hand
side of (\ref{conj}) is
\begin{align*}
  & lhs(w,z_1,z_2) = \\ & {1\over{q^{-1}w - z_1}}\left( {{z_1 -
        z_2}\over{(q^2 z_1 - z_2)(z_1 - z_2)}} - (q+q^{-1}){{z_1 -
        z_2}\over{(q^2 z_2 - z_1)(q^{-1}z_2 - q z_1)}} \right. \\ &
  \left.  - {{z_1 - z_2}\over{(q^2 z_2 - z_1)(z_1 - z_2)}}\right) \\ &
  + {1\over{q^{-1}z_1 - w}} \left( {{z_1 - z_2}\over{(q^2 z_1 -
        z_2)(q^{-1} z_2 - q^{-1}z_1)}} + (q+q^{-1}) {{z_1 -
        z_2}\over{(q^2 z_1 - z_2)(q^{-2} z_1 - z_2)}} \right. \\ &
  \left.  - {{z_1 - z_2}\over{(q^2 z_2 - z_1)(q^{-1}z_2 - q^{-1}z_1)}}
  \right) \\ & + (z_1 \leftrightarrow z_2) ,
\end{align*}
that is 
\begin{align*}
& {1\over{q^{-1}w - z_1}}\left( {1\over{q^2 z_1 - z_2}} 
- (q+q^{-1}){{z_1 - z_2}\over{(q^2 z_2 - z_1)(q^{-1}z_2 - q z_1)}}
- {1\over{q^2 z_2 - z_1}}\right) 
\\ & +
{1\over{q^{-1}z_1 - w}}
\left( - q {1\over{q^2 z_1 - z_2}} 
+ (q+q^{-1}) {{z_1 - z_2}\over{(q^2 z_1 - z_2)(q^{-2} z_1 - z_2)}}
+ q {1\over{q^2 z_2 - z_1}} \right) 
\\ & + (z_1 \leftrightarrow z_2) . 
\end{align*}

Now 
$$ {{z_1 - z_2}\over{(q^2 z_2 - z_1)(q^{-1}z_2 - q z_1)}} =
- {q\over{q^2 + 1}} ( {1\over{q^2 z_2 - z_1}} + {1\over{z_2 - q^2 z_1}}), 
$$
and 
$$
{{z_1 - z_2}\over{(q^2 z_1 - z_2)(q^{-2} z_1 - z_2)}} = 
{q^2\over{q^2 + 1}} ( {1\over{q^2 z_1 - z_2}} + {1\over{z_1 - q^2 z_2}}),
$$
so that $lhs(w,z_1,z_2)$ is equal to 
\begin{align*} 
  & {1\over{q^{-1}w - z_1}}\left( {1\over{q^2 z_1 - z_2}} +
    {1\over{z_2 - q^2 z_1}}\right) + {1\over{q^{-1}z_1 - w}} q
  \left({1\over{q^2 z_2 - z_1}} + {1\over{z_1 - q^2 z_2}}\right) \\ &
  + (z_1 \leftrightarrow z_2) , 
\end{align*}
or
$$
{1\over{q^{-1}w - z_1}} \delta(q^2 z_1 , z_2) + q {1\over{q^{-1}z_1 - w}} 
\delta(q^2 z_2, z_1) 
  + (z_1 \leftrightarrow z_2) , 
$$
which is 
$$ q \delta(w, qz_1) \delta(q^2 z_1, z_2) + q
\delta(q^{-1}z_1,w)\delta(q^2z_2,z_1) 
$$ 
that is the right side of (\ref{conj}). 

In the case $m=2$, we expand the left side of (\ref{conj}) as 
\begin{equation} \label{step}
  \sum_{i=1}^3 {1\over{q^{-2}z_i - w}} f_i(z_1,z_2,z_3) + \sum_{i=1}^3
  {1\over{q^{-2}w - z_i}} g_i(z_1,z_2,z_3).
\end{equation}
Then $f_1(z_2,z_2,z_3)$ is expanded as 
$$ {1\over{q^2 z_1 - z_2}} h_1(z_2,z_3) + {1\over{q^{-4} z_1 - z_2}}
h_2(z_2,z_3) + {1\over{q^2 z_2 - z_1}} h_3(z_2,z_3) +
(z_2\leftrightarrow z_3). 
$$
We then compute 
$$
h_1(z_2,z_3) = 0 , \quad h_2(z_2,z_3) = q^{-2}\delta(z_3,q^2 z_2), \quad
h_3(z_2,z_3) = q^2 \delta(z_2,q^2 z_3), 
$$
so that 
\begin{align*}
  {1\over{q^{-2}z_1 - w}} f_1(z_1,z_2,z_3) & = {1\over{q^{-2}z_1 - w}}
  \left( {1\over{q^{-4}z_1 - z_3}} q^{-2} \delta(z_2,q^2 z_3) \right.
      \\ & \left. + {1\over{q^2 z_2 - z_1}} q^2 \delta(z_2,q^2 z_3) +
    (z_2 \leftrightarrow z_3) \right) \\ & = {1\over{q^{-2}z_1 - w}}
  q^2 \delta(z_1,q^2 z_2)\delta(z_2,q^2 z_3) + (z_2 \leftrightarrow
  z_3).
\end{align*}
It follows that the terms associated with $i=1$ in (\ref{step}) are equal to 
$$ q^2 \delta(w,q^{-2}z_1)\delta(z_1,q^2 z_2)\delta(z_2,q^2 z_3) +
(z_2 \leftrightarrow z_3).
$$ Adding up the contributions of all $i$ we find that the left side
of (\ref{conj}) is equal to its right side.  \hfill \qed \medskip

\begin{remark} 
  It is natural to expect that identity (\ref{conj}) is valid for any
  $m$. For $m=3$, this would imply the statement of Thm. \ref{main}
  also for $\bar\G$ of type $G_2$.  Identity (\ref{conj}) implies the
  combinatorial identity (6.1) of \cite{Jing}.
\end{remark}

\section{Vanishing properties of the correlation functions}
\label{sect:proof:3}

Let us now show that identity (\ref{conj}) imply Thm.  \ref{main}, 3).

This statement is nonempty only if $N_{\al} \geq 1 - a_{\al\beta}$ and
$N_{\beta} \geq 1$.  If $a_{\al\beta} = 0$, the commutation of
$e_\al(z^{(\al)}_i)$ and $e_\beta(z^{(\beta)}_j)$ implies that
$z^{(\al)}_i - z^{(\beta)}_j$ divides $A(z^{(\al)}_i)$.

Assume that $a_{\al\beta}$ is equal to $-1$ or $-2$. Set $z_i =
z^{(\al)}_i$, for $i=1,\ldots, 1- a_{\al\beta}$ and $w =
z^{(\beta)}_1$.  Denote as $\prod_{i=1}^{N'}e_{\al_i}(z'_i)$ the
product
$$
\prod_{i>1-a_{\al\beta}} e_\al(z^{(\al)}_i) \prod_{i>1}
e_\beta(z^{(\beta)}_i) \prod_{\gamma\neq\al,\beta} \prod_{i\in
  I_\gamma}e_{\gamma}(z^{(\gamma)}_i)
$$
and set 
$$ v' = \prod_{i=1}^{N'}e_{\al_i}(z'_i) v
$$
and
$$ P = \prod_{i=1...k,j=1...N'}{{z_i - z'_j}\over{q^{(\al,\al_i)}z_i -
    z'_j}} \prod_{i=1}^{N'} {{w-z'_i}\over{q^{(\beta,\al_i)}w-z'_i}}
\prod_{1\leq i < j \leq N'} {{z'_i - z'_j}\over{q^{(\al_i,\al_j)}z'_i
    - z'_j}}, $$ where we denote by $(\al,\beta) = d_\al a_{\al\beta}$
the scalar product of two simple roots. $P$ belongs to
$\CC[[z_i,w]]((z'_1)) \cdots ((z'_{N'}))$.

Then we have for any $k = 0,\ldots, 1-a_{\al\beta}$,
\begin{align*}
  & \langle \xi, e_\al(z_1)\cdots e_\al(z_k) e_{\beta}(w)
  e_\al(z_{k+1}) \cdots e_\al(z_{1-a_{\al\beta}}) v' \rangle \\ & =
  (-1)^k\prod_{i=1}^{1-a_{\al\beta}} (w-z_i) \cdot P \cdot {
    {A(z^{(\al)}_i)} \over {\prod_{i=1}^k(q^{(\al,\beta)}z_i - w)
      \prod_{i=k+1}^{1-a_{\al\beta}}}(q^{(\al,\beta)}w - z_i)}
  \prod_{1\leq i<j\leq 1-a_{\al\beta}}{{z_i - z_j}\over{q_\al}^2 z_i -
    z_j}
\end{align*} 
[equality in
$\CC((z_1))\cdots((z_k))((w))((z_{k+1}))\cdots((z_{1-a_{\al\beta}}))
((z'_1))\cdots((z'_{N'}))$].

The quantum Serre relation (\ref{q:Serre}) implies that $\langle \xi,$
left side of relation (\ref{q:Serre}) $v'\rangle$ is equal to zero. 
(\ref{conj}) implies that this relation is written as 
$$ P \cdot \left(
  \on{Sym}_{z_1,\ldots,z_{1-a_{\al\beta}}}\delta(w,q_\al^{-a_{\al\beta}}z_1)
  \delta(z_1,q_\al^2 z_2) \cdots \delta(z_{-a_{\al\beta}},q_\al^2
  z_{1-a_{\al\beta}})\right) A(z^{(\al)}_i) = 0.
$$ This implies that the product of the last two terms is itself equal
to zero. But the product of $A(z^{(\al)}_i)$ with each delta-function
is the product of this delta-function and of a function depending only
on the $z'_i$ and $w$. Since the delta-functions are linearly
independent, the evaluation of $A(z^{(\al)}_i)$ on each variety
$z^{(\al)}_{\si(1)} = q_\al^2 z^{(\al)}_{\si(2)} = \cdots =
q_\al^{-2a_{\al\beta}} z^{(\al)}_{\si(1-a_{\al\beta})} =
q_\al^{-a_{\al\beta}} z_{1}^{(\beta)}$ is zero. Since $A$ is symmetric
in the group of variables $z^{(\al)}_i$, this is equivalent to the
conditions of Thm. \ref{main}, 3).
This completes the proof of Thm. \ref{main}.

\begin{remark}
  Conditions (\ref{vanishing}) were obtained in \cite{DM} in the case
  where $\bar\G = \SL_n$ and $V$ is integrable, using products of
  Frenkel-Kac realizations.  In that case, the correlation functions
  (\ref{def:corr}) have other functional properties that were studied
  in that paper.
\end{remark}

\section{Application to functional description of $U_q\N_+$}
\label{sect:fun:desc}

In this section, we prove Thm. \ref{thm:shuffle}.

Define dual Hopf algebras $(U_q\wt \B_\pm,\Delta_\pm)$ as follows.
$U_q\B_+$ has generators $e_\al[n]$, $n\in \ZZ$ and $K^+_\al[n],n\ge
0$, $K^+_\al[0]^{-1}$, generating series $e_\al(z) = \sum_{n\in\ZZ}
e_\al[n]z^{-n}$, $K^+_\al(z) = \sum_{n\ge 0} K^+_\al[n]z^{-n}$,
relations (\ref{vertex}) and
$$
[K^+_\al[n], K^+_\beta[m]] = 0, 
$$
and
$$ K^+_\al(z) e_\beta(w) K^+_\al(z)^{-1} =
{{z-q^{(\al,\beta)}w}\over{q^{(\al,\beta)} z-w}} e_\beta(w);
$$
(expansion for $w<<z$); 
the coproduct is defined by 
$$ \Delta_+(K^+_\al(z)) = K^+_\al(z) \otimes K^+_\al(z), \quad
\Delta_+(e_\al(z)) = e_\al(z) \otimes K^+_\al(z) + 1 \otimes e_\al(z).
$$

$U_q\B_-$ has generators $f_\al[n]$, $n\in \ZZ$ and $K^-_\al[n],n\leq
0$, $K^-_\al[0]^{-1}$, generating series $f_\al(z) = \sum_{n\in\ZZ}
f_\al[n]z^{-n}$, $K^-_\al(z) = \sum_{n\ge 0} K^-_\al[-n]z^{n}$,
relations (\ref{vertex}) with $q$ replaced by $q^{-1}$ between the
$f_\al(z)$ and
$$
[K^-_\al[n], K^-_\beta[m]] = 0, 
$$
and
$$ K^-_\al(z) f_\beta(w) K^-_\al(z)^{-1} =
{{z-q^{(\al,\beta)}w}\over{q^{(\al,\beta)} z-w}} f_\beta(w);
$$
(expansion for $z<<w$); 
the coproduct is defined by 
$$ \Delta_-(K^-_\al(z)) = K^-_\al(z) \otimes K^-_\al(z), \quad
\Delta_-(f_\al(z)) = f_\al(z) \otimes 1 + K^-_\al(z)^{-1}\otimes
f_\al(z).
$$ So $U_q\wt\B_\pm$ are the usual opposite Hopf subalgebras of the
new realizations algebras (\cite{D}), where the quantum Serre
conditions are not imposed. The following result can be viewed as an
infinite-dimensional analogue of results in \cite{L,Ro}.

\begin{prop} \label{serre:to:0}
  We have a Hopf algebra pairing between $(U_q\wt\B_+,\Delta_+)$ and
  $(U_q\wt\B_-,\Delta'_-)$, defined by
  $$\langle K^+_\al(z),K^-_\beta(w) \rangle =
  {{z-q^{(\al,\beta)}w}\over{q^{(\al,\beta)}z-w}}, \quad \langle
  e_\al[n],f_\beta[m] \rangle = \delta_{\al\beta}\delta_{n+m,0}$$
  (expansion for $w<<z$).  The ideals defined by the quantum Serre
  relations are contained in the radicals of this pairing.
\end{prop}

{\em Proof.} The verification of the first statement is standard.  To
show the last statement, let us compute the pairing of the Serre
relation (\ref{q:Serre}) with any element of $U_q\wt\B_-$.
$U_q\wt\B_-$ has a gradation by the root lattice of $\bar\G$, defined
by deg$(f_\al[n]) = \al$ and deg$K_\al[-n] = 0$. Then the pairing of
the left side of (\ref{q:Serre}) can be nontrivial only against an
element of degree $(1-a_{\al\beta})\al + \beta$.  Translating the
Cartan modes $K^-_\al[n]$ to the left of the $f_{\al,\beta}[\phi]$, we
have to compute for any $k$, and $\phi_i,\psi$ in $\CC[z,z^{-1}]$,
\begin{equation} \label{oz} \langle \on{left\ side\ of\ }(\ref{q:Serre}),
  f_\al[\phi_1]\cdots f_\al[\phi_k]f_\beta[\psi]
  f_\al[\phi_{k+1}]\cdots f_\al[\phi_{1-a_{\al\beta}}] \rangle .
\end{equation} 
Denote by $L_m(z_1,\cdots,z_{m+1},w)$ the left side of identity
(\ref{conj}), with $q$ replaced by $q^{-1}$, viewed as a rational
function and expanded for $z_1>>z_2>>\cdots>>z_{m+1}>>w$.  We find
that (\ref{oz}) is equal to
\begin{align*} 
  & \on{Sym}_{z_1\cdots z_{1-a_{\al\beta}}} \left( \phi_1(z_1)\cdots
    \phi_{1-a_{\al\beta}}(z_{1-a_{\al\beta}}) \psi(w) \prod_{i=1}^k
    (q^{(\al,\beta)}z_i - w)\prod_{i=k+1}^{1-a_{\al\beta}}
    (q^{(\al,\beta)}w - z_i) \right. \\ & \left. \prod_{j>i}(z_j -
    q_\al^{-2}z_i) L_{-a_{\al\beta}}(z_1,\ldots,z_{1-a_{\al\beta}},w)
  \right).
\end{align*}
On the other hand, from \cite{Jing} follows that
$L_m(z_1,\cdots,z_{m+1},w)$ is identically zero (as a rational
function, and therefore as a formal series in
$\CC((z_1))\cdots((z_{m+1}))((w))$.  It follows that (\ref{oz}) is
zero. In the same way, one shows that the quantum Serre relations of
$U_q\wt\B_-$ are in the radical of the pairing.  \hfill \qed \medskip

\begin{cor} 
  The coproducts $\Delta_\pm$ induce Hopf algebra structures on the
  quotients $U_q\B_\pm$ of $U_q\wt\B_\pm$ by the ideals generated by
  the quantum Serre relations. 
\end{cor}

Then

\begin{prop}
  The direct sum $\wt\pi$ of the maps $\wt\pi_\nn$ from $U_q\wt\B_+$
  to $\CC[[z^{(\al)}_i,z^{(\al)-1}_i]]$, defined by $\wt\pi_{\nn}(x) =
  \langle x , \prod_{\al=1}^r \prod_{i=1}^{n_r} f_\al(z^{(\al)}_i)
  \rangle$ induces a linear map $\pi$ from $U_q\B_+$ to
  $\overline{Sh}$; if $\bar\G$ is not of type $G_2$, the image of
  $\pi$ is contained in $Sh$.
\end{prop}

{\em Proof.} The Hopf pairing rules and the commutation relations
imply that $\wt\pi_\nn$ defines a linear map from $U_q\wt \B_+$ to
$\CC((z^{(1)}_1)) \cdots ((z^{(r)}_{n_r}))$.  Crossed vertex relations
allow to apply the reasoning of sect. \ref{sect:std} to prove exchange
relations, which imply that the image of $\wt\pi_\nn$ is contained in
subspace of elements of the form (\ref{can}), with $A$ in
$\CC[[z^{(\al)}_i]][z^{(\al)-1}_i]$.  A degree argument can be used as
in sect. \ref{sect:std} to show that $A$ belongs to $\CC[z^{(\al)}_i,
z^{(\al)-1}_i]$.  Therefore $\wt\pi$ maps $U_q\wt\B_+$ to
$\overline{Sh}$. Since $\wt\pi$ sends the radical of $\langle,\rangle$
to zero, it induces a map $\pi$ from $U_q\B_+$ to $\overline{Sh}$.

If $\bar\G$ is not of type $G_2$, we can then follow the reasoning of
sect. \ref{sect:proof:3} to show that $A$ satisfies (\ref{vanishing}),
which shows that the image of $\wt\pi$ is contained in $Sh$. \hfill
\qed \medskip

Define the shuffle product on $\overline{Sh}$ by the following rule.
Let $f$ and $g$ belong to in $\overline{Sh}_\nn$ and
$\overline{Sh}_\mm$. Set $z_1 = z^{(1)}_1, z_2 = z^{(1)}_2$, etc.,
$z_{n_1 +m_1+ \cdots + n_{\al - 1} +m_{\al-1}+ i} = z^{(\al)}_i$ for
$i = 1 , \ldots, n_\al+m_\al$, $z_{N+M} = z^{(r)}_{n_r+m_r}$; we set
$N = \sum_{i=1}^r n_i$, $M = \sum_{i=1}^r m_i$.  We associate
``colors'' to the variables $z_i$ by the rule $\al(i) = \al$ if $i =
n_1 + m_1 + \cdots + n_{\al-1} + m_{\al-1} + i$, $i =
1,\cdots,n_\al+m_\al$.  Let $M_{\nn,\mm}$ be the set of bijective maps
from $\{1,\ldots,N+M\}$ to itself, such that for each $\al$, we have
\begin{align*} 
\{ \sum_{i=1}^{\al-1} (n_i + m_i) + 1, \ldots, \sum_{i=1}^{\al} (n_i +
m_i) \} & = \sigma( \{\sum_{i=1}^{\al-1} n_i + 1, \ldots,
\sum_{i=1}^{\al} n_i \} ) \\ & \cup \sigma(\{N+\sum_{i=1}^{\al-1} m_i +
1, \ldots, N+\sum_{i=1}^{\al} m_i\})
\end{align*} and $\si(i) < \si(j)$ if $i<j$ and $i,j$ both belong to some $
\{\sum_{i=1}^{\al-1} n_i + 1, \ldots, \sum_{i=1}^{\al} n_i \}$ or
$\{N+\sum_{i=1}^{\al-1} m_i + 1, \ldots, N+\sum_{i=1}^{\al} m_i\}$.
Then the product $f*g$ is defined in $\overline{Sh}_{\nn+\mm}$ as
\begin{align} \label{sh:pdt}
  (f*g)(z_1,\ldots,z_{N+M}) & = \sum_{\sigma\in M_{\nn,\mm}}
  \varepsilon(\sigma) \prod_{i<j, \si^{-1}(i)>\si^{-1}(j)}
  {{q^{(\al(i),\al(j))}z_i - z_j}\over{z_i - q^{(\al(i),\al(j))} z_j}}
  \\ & \nonumber
  f(z_{\si(1)},\ldots,z_{\si(N)})g(z_{\si(N+1)},\ldots,z_{\si(N+M)}).
\end{align}
Note that each summand in the right side belongs to $\CC((z_1))\cdots
((z_{N+M}))$, because the prefactor cancels poles in $f$ or $g$.

\begin{remark} 
  $M_{\nn,\mm}$ is isomorphic to the product $\prod_{i=1}^r
  sh_{n_i,m_i}$, where $sh_{n,m}$ is the set of shuffle maps of the
  pair of sequences $((1,\ldots,n),(n+1,\ldots,n+m))$.
\end{remark}

\begin{remark} The product (\ref{sh:pdt}) is a transcription of the shuffle 
  algebra structure $(FO,\wt *)$ defined in \cite{FO}. The isomorphism
  $i$ from $Sh$ to $FO$ is given by the formula
  $$ i(f)(z_1,\ldots,z_N) = \prod_{i< j}{{z_i -
      q^{(\al(i),\al(j))}z_j} \over{z_i - z_j}} f(z_1,\ldots,z_N).
  $$ 
\end{remark}

\begin{prop}
  For $\bar\G$ arbitrary, $\pi$ is a morphism of algebras from $U_q\N_+$
  to $\overline{Sh}$, endowed with the shuffle product (\ref{sh:pdt}).
\end{prop}

{\em Proof.} It suffices to see that $\wt \pi$ is an algebra morphism.
It is clear that $\wt\pi$ maps $e_\al[n]$ to the element
$(z^{(\al)})^n$ of $\overline{Sh}_{\delta_\al}$. So we should check
that $\wt \pi$ maps the crossed vertex relations to zero. This is an
easy computation, relying on the equalities
$$
(z^{(\al)})^n * (z^{(\beta)})^m = (z^{(\al)})^n(z^{(\beta)})^m
$$
and
$$ (z^{(\beta)})^m * (z^{(\al)})^n = {{q^{(\al,\beta)}z^{(\al)} -
      z^{(\beta)}} \over{{z^{(\al)} - q^{(\al,\beta)}
      z^{(\beta)}}}}(z^{(\al)})^n(z^{(\beta)})^m, 
$$
for $\al < \beta$, and 
$$ (z^{(\al)})^n * (z^{(\al)})^m = (z_1^{(\al)})^n (z_2^{(\al)})^m +
{{q_\al^2 z_1 - z_2 }\over {z_1 - q_\al^2 z_2 }} (z_2^{(\al)})^n
(z_1^{(\al)})^m .
$$
\hfill \qed\medskip 

\begin{remark}
  Conditions (\ref{vanishing}) appeared in \cite{FOJM} in an elliptic
  situation and the $\SL_n$ case as a substitute to quantum Serre
  relations. 
\end{remark}

\end{document}